\documentstyle{article}\begin{document}\centerline{\bf An integral representation for the product of two parabolic cylinder functions }\centerline{\bf having unrelated arguments.}\vskip 1in

\centerline{M.L. Glasser}

\centerline{Department of Physics, Clarkson University, Potsdam,\\
New York 13699-5820, USA}

\centerline{Donostia International
Physics Center, P. Manuel de Lardizabal 4, \\ E-20018 San
Sebasti\'an, Spain}
\vskip .6in
\centerline{\bf Abstract}\vskip .1in

\begin{quote}
 An integral representation is provided for the parabolic cylinder function product $D_{\mu}(x)D_{\mu}(-y)$ where $Re\,\mu<0$  and $x>y$ are unrelated. A few simple consequences are given in the form of hyperbolic integrals and a sum rule.
 
 \end{quote}


\vskip .4in

\noindent
{\bf Class}

33C15, C3345, C3324

\newpage

\noindent
{\bf 1. Introduction}

A search of the literature and tables has shown that the number of known integral representations for a product of two parabolic cylinder functions is limited. Only two are available in the standard reference[1]
$$D_{\nu}(x)D_{-\nu-1}(x)=-\frac{1}{\sqrt{\pi}}\int_0^{\infty}\coth^{\nu+1/2}\frac{t}{2}\frac{1}{\sqrt{\sinh t}}\sin\frac{x^2\sinh t+\pi\nu}{2}dt$$

$$D_{\nu}(ze^{i\pi/4})D_{\nu}(ze^{-i\pi/4})=\frac{1}{\Gamma(-\nu)}\int_0^{\infty}\coth^{\nu} t\exp\left[-\frac{z^2}{2}\sinh 2t\right]\frac{dt}{\sinh t}.$$
and  three more are given in the Wolfram web-site[2]
$$D_{\nu}(x)D_{-\nu-1}(x)=2\int_0^{\infty}e^{-xt}\cos\left(xt-\frac{\pi\nu}{2}\right)J_{\nu+1/2}(t^2)dt$$
$$D_{\nu}(ze^{i\pi/4})D_{\nu}(ze^{-i\pi/4})=\frac{\sqrt{\pi}}{\Gamma(-\nu)}\int_0^{\infty}e^{-zt}J_{-\nu-1/2}(\frac{t^2}{2})dt$$

$$D_{\nu}(ze^{i\pi/4})D_{\nu}(ze^{-i\pi/4})=\frac{2\sqrt{2}}{\sqrt{\pi}\Gamma(-\nu)}\int_0^{\infty}e^{-zt}\cos\left(zt-\frac{\nu\pi}{2}\right)K_{\nu+1/2}(t^2)dt,$$
where, in terms of the Whittaker function,
$$D_{\nu}(z)=2^{(2\nu+1)/4}z^{-1/2}W_{{2\nu+1)/4},-1/4}(z^2/2).$$
Even these five are  inter-related by contour manipulation. A fairly  recent report by C. Malyshev[3] contains numerous references to the origins of these representations, to which we refer the interested reader. Conspicuously missing among these entries are representations where the two arguments on the left hand side are not linearly related. The aim of the present note is to fill this gap, at least partially since representations where the indices are unrelated is still open.

\vskip .2in
\noindent
{\bf 2. Derivation}

\vskip .1in
We begin by examining the Laplace transform
$$I=\int_0^{\infty}t^{\nu-1}(1+t)^{-\nu-1/2}e^{-at}e^{b\sqrt{t(t+1)}}dt,\quad \nu>0,\quad Re[a-b]>0.\eqno(1)$$
By introducing the substitution $t=u^2/(1-u^2)$ one finds
$$I=2\int_0^1u^{2\nu-1}(1-u^2)^{\nu-3/2}\exp\left\{\frac{bu-au^2}{1-u^2}\right\}du.\eqno(2)$$
Next, we introduce the quantities $X=\sqrt{a+\sqrt{a^2-b^2}}/\sqrt{2}$, $Y=\sqrt{a-\sqrt{a^2-b^2}}/\sqrt{2}$ and recall Mehler's formula [4] 
$$\exp\left[\frac{2XYu-(X^2+Y^2)u^2}{1-u^2}\right]=\sqrt{1-u^2}\sum_{n=0}^{\infty}\frac{H_n(X)H_n(Y)}{2^nn!}u^n\eqno(3)$$
to find, after performing the elementary $u$-integration,
$$I=2\sum_{n=0}^{\infty}\frac{H_n(X)H_n(Y)}{2^nn!(2\nu+n)} .\eqno(4)$$

\vskip .2in
   Now let us consider the Sturm-Liouville problem on the real line
   $$y''+(\lambda-x^2)y(x)=0,\quad y(\pm\infty)=0\eqno(5)$$
   for which the  normalized eigenfunctions are 
   $$y_n(x)=\frac{1}{\sqrt{2^nn!\sqrt{\pi}}}e^{-\frac{1}{2}x^2}H_n(x), \quad \lambda_n=2n+1,\quad n=0,1,2,\cdots.\eqno(6)$$
   Then the solution to the Green function equation
   $$Y''(x,x',\lambda)+(\lambda-x^2)Y(x,x',\lambda)=\delta(x-x'), \quad Y(\pm\infty)=0\eqno(7)$$
   is
   $$Y(x,x',\lambda)=\frac{1}{\sqrt{\pi}}e^{-\frac{1}{2}(x^2+x'^2)}\sum_{n=0}^{\infty}\frac{H_n(x)H_n(x')}{2^nn!(\lambda_n-\lambda)}.\eqno(8)$$
   However, the problem (5) has been treated by Titchmarsh[4] who has shown that, for $x> x'$,
   $$Y(x,x',\lambda)=\frac{1}{2\sqrt{\pi}}\Gamma\left(\frac{1-\lambda}{2}\right)D_{\frac{\lambda-1}{2}}(x\sqrt{2})D_{\frac{\lambda-1}{2}}(-x'\sqrt{2}).\eqno(9)$$
   Finally, from (4), (8) and (9) we conclude that  for $x>y>0$ and $\nu>0$,
   $$D_{-\nu}( x)D_{-\nu}(-y)$$
   $$=\frac{e^{-\frac{1}{4}(x^2+y^2)}}{2\Gamma(\nu)}\int_0^{\infty}t^{\nu/2-1}(t+1)^{-(\nu+1)/2}e^{-\frac{1}{2}(x^2+y^2)t+xy\sqrt{t(t+1)}}dt.\eqno(10)$$
   \vskip .2in
   \noindent
   {\bf 3.Discussion}\vskip .1in
   
   Equation (10) is our principal result and can be written
   $$\int_0^{\infty}t^{\nu/2-1}(1+t)^{-(\nu+1)/2}e^{-at} e^{b\sqrt{t(t+1)}}dt$$
   $$=2e^{a/2}\Gamma(\nu)D_{-\nu}\left(\sqrt{a+\sqrt{a^2-b^2}}\right)D_{-\nu}\left(-\sqrt{a-\sqrt{a^2-b^2}}\right),\eqno(11)$$
   for $Re(a)>Re(b)>0$. Similarly, it can be proven that
   $$\int_0^{\infty}t^{\nu/2-1}(1+t)^{-(\nu+1)/2}e^{-at} e^{-b\sqrt{t(t+1)}}dt$$
   $$=2e^{a/2}\Gamma(\nu)D_{-\nu}\left(\sqrt{a+\sqrt{a^2-b^2}}\right)D_{-\nu}\left(\sqrt{a-\sqrt{a^2-b^2}}\right),\eqno(12)$$
   for $ Re(a+b)>0$..

   First of all, (10) cannot be related to the  identities listed in the introduction, for the right hand side diverges when $x=y$. In the forms (11) and  (12), since both sides can be differentiated  any number of times with respect to the parameters $a$, $b$ or $\nu$, it is the gateway to a previously unknown class of Laplace Transforms. Furthermore, through (4) we have summed exactly a new Hermite series from which  further ones can be obtained. To conclude, we append a short list of special cases.
   
     For $\nu=1$ , $D_{-1}(z)=\sqrt{\pi/2}\exp(z^2/4){\rm Erfc}(z/\sqrt{2})$ and (12) gives rise  to the interesting hyperbolic integrations
   $$\int_0^{\infty}{\rm sech}\theta\, e^{-\alpha^2\sinh\theta\sinh(\theta+\phi)}d\theta$$
   $$=\frac{\pi}{2}e^{\alpha^2\cosh\phi} {\rm Erfc}(\alpha\sinh(\phi/2))\rm{Erfc}(\alpha\cosh(\phi/2)).\eqno(13a)$$
   $$\int_0^{\infty}\sinh\theta\, e^{-\alpha^2\sinh\theta\sinh(\theta+\phi)}d\theta$$
   $$=\frac{\sqrt{\pi}}{2\alpha}\left[e^{\alpha^2\cosh^2(\phi/2)}\cosh(\phi/2){\rm Erfc}(\alpha\cosh(\phi/2))-e^{\alpha^2\sinh^2(\phi/2)}\,\sinh(\phi/2){\rm  Erfc}(\alpha\sinh(\phi/2))\right]\eqno(13b)$$
   
   In the case $\nu=1/2$, where $D_{-1/2}(z)=\sqrt{z/2\pi}K_{1/4}(z^2/4)$, (12) can be manipulated to give
   $$\int_0^{\infty}\frac{d\theta}{\sqrt{\sinh\theta}}e^{-a\cosh(\theta+\phi)}=\sqrt{\frac{a\sinh\phi}{\pi}} K_{\frac{1}{4}}(a\cosh^2(\phi/2))K_{\frac{1}{4}}(a\sinh^2(\phi/2)).\eqno(14)$$
   
   Finally, From (3) we obtain the new sum rule 
   $$\sum_{n=0}^{\infty}\frac{D_n(x)D_n(y)}{n!(n +\nu)}=\Gamma(\nu)D_{-\nu}(x)D_{-\nu}(-y),\mbox{  where  } x>y,\eqno(15)$$
   from which others may be derived.

   \newpage
   \noindent
   {\bf Acknowledgement}\vskip .1in
   
   The author is grateful for hospitality of the DIPC where this work was carried out.
   
   \centerline{\bf References}\vskip .1in
   \noindent
   [1]  I.S. Gradshteyn and I. Ryzhik, {\it Tables of Integrals, Series and Products}, Ed. A. Jeffrey and D. Zwillinger, Academic Press, N.Y. (2007).
   
   \noindent
   [2] Wolfram Website: http://functions.wolfram.com/hypergeometricfunctions/paraboliccylinderD/07/01/
   
   \noindent
   [3] C. Malyshev,  arXiv:math/01/06142v2[math.CA]4Dec2001
   
   \noindent
   [4]  Weisstein, Eric W. "Mehler's Hermite Polynomial Formula." From MathWorld--A Wolfram Web Resource. http://mathworld.wolfram.com/MehlersHermitePolynomialFormula.html 
   
   \noindent
   [4] E.C. Titchmarch, {\it Eigenfunction Expansions, Vol.1}, Oxford, London (1961) p.74.

\end{document}